\documentclass [12pt,a4paper]{article}
\usepackage{amssymb, theorem}
\newtheorem{thm}{Theorem}
\newtheorem{cor}{Corollary}

\newtheorem{lemma}{Lemma}
\theorembodyfont{\rm}
\newtheorem{pl}{Example}

\def\proof{{\it Proof: }}

\def\Imag{\mathrm{Im}\,}

\def\qed{\nobreak\hfill $\square$}
\def\<{\langle}
\def\>{\rangle}

\def\bL{{\bf L}}

\def\bM{{\bf M}}

\def\bbbr{{\mathbb R}}

\def\Diag{\mbox{Diag}\,}

\def\Tr{\mathrm{Tr}\,}

\def\J{{\mathbb J}}

\def\bL{{\mathbb L}}
\def\bR{{ \mathbb R}}

\begin{document}
%\rightline{\today }
\ \vskip 1cm 
\centerline{\LARGE {\bf Completely positive mappings}}
\bigskip
\centerline{\LARGE {\bf and mean matrices}}
\bigskip
\bigskip
%\centerline{\bf Not for circulation}
\bigskip
\bigskip
\centerline{{\bf \'Ad\'am Besenyei\footnote{E-mail: badam@cs.elte.hu.} and
D\'enes Petz\footnote{E-mail: petz@math.bme.hu.}}}
\bigskip
%Partially supported by the Hungarian Research Grant OTKA  T068258.}}
%\begin{center}
\centerline{Department of Applied Analysis}
\centerline{ELTE, H-1117 Budapest, P\'azm\'any s\'et\'any 1/c, Budapest, Hungary}
\medskip
\centerline{and}
\medskip
\centerline{Department of Mathematical Analysis}
\centerline{ BME, H-1111 M\H uegyetem rakpart 3-9, Budapest, Hungary}
%\centerline{Alfr\'ed R\'enyi Institute of Mathematics, H-1364 Budapest,
%POB 127, Hungary}
\bigskip
\begin{abstract}
Some functions $f:\bbbr^+ \to \bbbr^+$ induce mean of positive numbers and the
matrix monotonicity gives a possibility for means of positive definite matrices.
Moreover, such a function $f$ can define linear mapping $\beta_f:\bM_n \to \bM_n$
on matrices (which is basical in the constructions of monotone metrics). The present
subject is to check the complete positivity of $\beta_f$ in the case of a few concrete 
functions $f$. This problem has been motivated by applications in quantum information.


\end{abstract}

\bigskip
%%%%%%%%%%%%%%%%%%%%%
\section{Introduction}

The matrix monotone function $f:\bbbr^+ \to \bbbr^+$ will be called {\it 
standard} if $f(1)=1$ and $t f(t^{-1}) = f(t)$. Standard functions are used
to define (symmetric) matrix means:
$$
M_f(A,B)=A^{1/2}f(A^{-1/2}B A^{-1/2})A^{1/2},
$$
see \cite{K-A}. For numbers $m_f(x,y)= x f(y/x)$.

It is well-known that if $f:\bbbr^+ \to \bbbr^+$ is a standard matrix 
monotone function, then
$$
\frac{2x}{x+1} \le f(x) \le \frac{x+1}{2}.
$$
For example, 
$$
\frac{2x}{x+1} \le \sqrt{x} \le  \frac{x-1}{\log x} \le \frac{x+1}{2}
$$
they correspond to the harmonic, geometric, logarithmic and arithmetic
mean. The matrix means have application in quantum theory and this paper is
also motivated by that, see \cite{HP, PD22, PD143, TRus}.

Assume that a standard matrix monotone function $f$ is given. Let $\lambda_1, 
\lambda_2, \dots, \lambda_n$ be strictly positive numbers. The positivity of 
the matrix $X \in \bM_n$  defined as
\begin{equation}
X_{ij}=m_f(\lambda_i,\lambda_j)
\end{equation}
is an interesting question. We call $X$ {\it mean matrix}. Positivity of the mean matrix 
for all possibilities is equivalent to the positive definiteness of the kernel
$$
F_f(x,y):=m_f(x,y).
$$ 
(A stronger property than positivity is the so-called {\it infinite divisibility} 
\cite{Kos}, it is not studied here, but some results are used.)

The choice $\lambda_1=1$ and $\lambda_2=x$ shows that
$$
f(x) \le \sqrt{x}
$$
is a necessary condition for the positivity of the mean matrix, in other words 
$m_f$ should be smaller than the geometric mean. If $f(x) \ge \sqrt{x}$, then 
the matrix
\begin{equation}\label{E:T}
T_{ij}=\frac{1}{m_f(\lambda_i, \lambda_j)}
\end{equation}
can be positive. The matrix (\ref{E:T})
was important in the paper \cite{PD2} for the characterization of monotone 
metrics, see also \cite{PD22, PD143}. It will be shown that $T$ is positive
if and only if the linear mapping $ A \mapsto A \circ C$ is completely positive.
(Here $A \circ C$ is a notation for the Hadamard product.) 

The subject of the paper is the study of the existence and description of this 
kind of completely positive mappings which are induced by a standard matrix monotone 
function. Examples of good matrix means are presented.

%%%%%%%%%%%%%%%%%%%%%%%%%%%%%%%%%%%%%%%
\section{The positive operator  $\J_D^f$}

Let $D \in \bM_n $ be a positive definite matrix and $f$ be a standard matrix monotone
function. A linear operator $\J_D^f:  \bM_n \to  \bM_n $ can be defined.  If $D
=\Diag (\lambda_1,\dots, \lambda_n)$, then
$$
(\J_D^f A)_{ij}=A_{ij} m(\lambda_i,\lambda_j) \qquad (A \in \bM_n).
$$
Since 
$$
\Tr A^* (\J_D^f A)=\sum_{ij}|A_{ij}|^2 m(\lambda_i,\lambda_j) \ge 0,
$$
the linear mapping $\J_D^f$ is positive with respect to the Hilbert-Schmidt inner product
(for any $f\ge 0$). Another definition is
\begin{equation}
\J_D^f=f(\bL_D\bR_D^{-1})\bR_D\,,
\end{equation}
where
$$
\bL_D(X)=D X \qquad\mbox{and}\qquad \bR_D(X)=XD\,.
$$
(The operator $\bL_D\bR_D^{-1}$ appeared in the modular theory of von Neumann
algebras.)

The inverse of this mapping is
$$
((\J_D^f)^{-1} B)_{ij}=B_{ij} \frac{1}{m(\lambda_i,\lambda_j)} \qquad (A \in \bM_n).
$$
and it appeared in \cite{PD2} to describe the abstract quantum Fisher information,
see also \cite{PD143}. The linear mappings $(\J_D^f)^{-1}: \bM_n \to \bM_n$ have
the monotonicity condition
\begin{equation}\label{E:Fmon}
\alpha^* (\J_{\alpha(D)}^f)^{-1}\alpha  \le (\J_D^f)^{-1}
\end{equation} 
for every completely positive trace preserving mapping $\alpha : \bM_n \to \bM_m$, if
$f$ is a standard matrix monotone function.

The linear transformation $(\J_D^f)^{-1}$ appeared also in the paper \cite{TRus} in a different
notation. There $\Omega^k_D$ is the same as  $(\J_D^f)^{-1}$ with $f=1/k$, see also 
\cite{PD22, PD143}. The complete positivity of the mapping $\beta:=(\J_D^f)^{-1}:
\bM_n \to \bM_n$ is a question in the paper \cite{TRus}. 

The subject of this paper is to find functions $f$ such that the mapping $\beta$ is monotone 
(in the sense of (\ref{E:Fmon}) ) and completely positive. The complete positivity is
equivalent to the positivity of a matrix, see the next lemma. 

If $D=\Diag (\lambda_1,\dots, \lambda_n)$, then
$$
\beta (A)_{ij}=A_{ij} \frac{1}{m(\lambda_i,\lambda_j)},
$$
where $m$ is the mean corresponding to the function $f$. In the notation (\ref{E:T}),
we have $\beta (A)=A \circ T$, it is a Hadamard product. 

\begin{lemma}
The linear mapping $\beta :\bM_n \to \bM_n, \beta(A)=A \circ T$ is completely positive
if and only if the matrices $T \in \bM_n$ defined in (\ref{E:T}) are positive.
\end{lemma}

\proof
If $\beta$ is completely positive, then $ A \circ T \ge 0$ for every positive $A$.
This implies the positivity of $T$.

The mapping $ \beta$ linearly depends on $T$. Therefore, it is enough to prove
the complete positivity when  $T_{ij}= \overline \lambda_i \lambda_j$. Then
$$
\beta(A)=\Diag(\lambda_1,\lambda_2,\dots,\lambda_n)^*A\Diag(\lambda_1,\lambda_2,
\dots,\lambda_n)
$$
and the complete positivity is clear. \qed

Hence $\beta$ is  completely positive if and only if $\beta$ is positive. The problem 
of \cite{TRus} is equivalent to the positivity of the matrix $T$.

%%%%%%%%%%%%%%%%%%%%%%%%%%%%%%%%%%%%%
\section{Completely positive mappings}

In this section we analyze the complete positivity of $\beta=\J_D^f$ for several 
matrix monotone functions $f$. The first three examples are very simple and actually 
they are particular cases of Example \ref{Pl:6}.

\begin{pl}
If $f(x)=\sqrt{x}$, then
$$
\< a, Ta\>=\sum_{ij} \frac{\overline a_i a_j}{m(\lambda_i, \lambda_j)}=
\overline{\sum_{i} a_i \lambda_i^{-1/2}}\sum_{j} a_j \lambda_j^{-1/2} \ge 0,
$$
so $T \ge 0$. In this case $\beta (A)= D^{-1/2}A D^{-1/2}$ and the complete positivity
is obvious. Moreover,  $\beta^{-1}$ is completely positive as well. (This is the only
example such that both $\beta$ and $\beta^{-1}$ are completely positive.) \qed
\end{pl}

\begin{pl}
If $f(x)=(1+x)/2$, the arithmetic mean, then $T$ is the so-called Cauchy matrix,
$$
T_{ij}=\frac{2}{\lambda_i+\lambda_j}=2 \int_0^\infty e^{s \lambda_i}e^{s \lambda_j}\, ds,
$$
which is positive. Therefore $\beta: A \mapsto A \circ T$ is completely positive. This
can be seen also from the formula
$$
\beta(A)=2 \int_0^\infty \exp (-sD)A\exp (-sD)\, ds.
$$
\qed
\end{pl}

\begin{pl}\label{P:3}
The {\it logarithmic mean} corresponds to the function $f(x)=(x-1)/\log x$.

Let $D=\Diag(\lambda_1, \lambda_2,\dots,\lambda_n)$ be positive definite. The mapping
$$
\beta: A \mapsto \int_0^\infty (D+t)^{-1}A (D+t)^{-1}\,dt
$$
is a positive mapping. Since $\beta(A)= T\circ A$ is a Hadamard product with
$$
T_{ij}=\frac{\log \lambda_i - \log \lambda_j}{\lambda_i-\lambda_j},
$$
the positivity of the mapping $\beta$ implies the positivity of $T$. Another proof
comes from the formula
$$
\frac{\log \lambda_i - \log \lambda_j}{\lambda_i-\lambda_j}=
\int_0^\infty\frac{1}{(s+ \lambda_i)(s+ \lambda_j)}\,ds.
$$
\qed
\end{pl}

\begin{pl}
Consider the mean 
$$
m(x,y):=\frac{1}{2}(x^ty^{1-t}+x^{1-t}y^t) \ge \sqrt{xy} \qquad (0<t<1)
$$
(which is sometimes called {\it Heinz mean}). Let $D=\Diag(\lambda_1, \lambda_2,\dots,
\lambda_n)$ be positive definite. 
The mapping
$$
\alpha: A \mapsto\frac{1}{2}( D^{t}AD^{1-t}+D^{1-t}AD^{t})
$$
has the form $\alpha (A)= X\circ A$ with
$$
X_{ij}=\frac{1}{2}(\lambda_i^t\lambda_j^{1-t}+\lambda_i^{1-t}\lambda_j^t).
$$
The inverse of the mapping is denoted by $\beta$, it is the Hadamard product with
$$
T_{ij}=\frac{2}{\lambda_i^t\lambda_j^{1-t}+\lambda_i^{1-t}\lambda_j^t}.
$$
$\beta$ is a positive mapping if and only $T\ge 0$.

To find the inverse of $\alpha$, we should solve the equation
$$
2A=D^{t}YD^{1-t}+D^{1-t}YD^{t},
$$
when $Y=\beta(A)$ is unknown. It has the form
$$
2D^{-t}AD^{-t}=YD^{1-2t}+D^{1-2t}Y
$$
which is a Sylvester equation. The solution is
$$
\beta(A)=Y= \int_0^\infty \exp(-s D^{1-2t})(2D^{-t}AD^{-t})\exp(-s D^{1-2t})\,ds.
$$
Therefore the mapping $\beta$ is positive and the matrix $T$ is positive as well. \qed
\end{pl}

The function
\begin{equation}\label{E:ft}
f_t(x)= 2^{2t -1}x^t (1+x)^{1-2t}
\end{equation}
is a kind of interpolation between the arithmetic mean  ($t=0$)
and the harmonic mean ($t=1$). This function appeared in the paper \cite{Han}
and it is proven there that it is a standard matrix monotone function.

\begin{thm}
If $t \in (0, 1/2)$, then
$$
f_t(x)= 2^{2t -1}x^t (1+x)^{1-2t} \ge \sqrt{x}
$$ 
and the matrix
$$
T_{ij}=\frac{1}{m_{f_t}(\lambda_i, \lambda_j)}=
\frac{2^{1-2t}}{(\lambda_i+\lambda_j)^{1-2t}}(\lambda_i\lambda_j)^{-t}
$$
is positive and the corresponding mapping $\beta$ is completely positive.  
\end{thm}

\proof
For $|x|<1$ and $1-2t=\alpha>0$ the binomial expansion yields
\[
(1-x)^{-\alpha}=\sum_{k=0}^\infty a_kx^k,
\]
where 
\[
a_k=(-1)^k {-\alpha \choose k}
=(-1)^k \frac{(-\alpha-1)(-\alpha-2)\cdot\dots\cdot(-\alpha-k+1)}{k!}>0.
\]
So that 
\begin{eqnarray*}
(\lambda_i+\lambda_j)^{-(1-2t)}&=&
\left(\left(\lambda_i+\frac12\right)\left(\lambda_j+\frac12\right)\left(1-
\frac{\left(\lambda_i-\frac12\right)\left(\lambda_j-\frac12\right)}{\left(\lambda_i
+\frac12\right)\left(\lambda_j+\frac12\right)}\right)\right)^{-(1-2t)}\cr &=&
\left(\lambda_i+\frac12\right)^{-(1-2t)}\left(\lambda_j+\frac12\right)^{-(1-2t)}\sum_{k=0}^\infty a_k\left(\frac{\left(\lambda_i-\frac12\right)\left(\lambda_j-\frac12\right)}{\left(\lambda_i+\frac12\right)\left(\lambda_j+\frac12\right)}\right)^k\cr&=&
\sum_{k=0}^\infty a_k\frac{\left(\lambda_i-\frac12\right)^k\left(\lambda_j-
\frac12\right)^k}{\left(\lambda_i+\frac12\right)^{k+(1-2t)}
\left(\lambda_j+\frac12\right)^{k+(1-2t)}}.
\end{eqnarray*}
Hence we have
\[
T_{ij}=2^{1-2t}\sum_{k=0}^\infty a_k\frac{\left(\lambda_i-\frac12\right)^k}
{\left(\lambda_i+\frac12\right)^{k+(1-2t)}\lambda_i^t}\,
\frac{\left(\lambda_j-\frac12\right)^k}{\left(\lambda_j+\frac12\right)^{k+(1-2t)}\lambda_j^t}
\]
and $T$ is the sum  of positive-semidefinite matrices of  rank one. \qed

If $t \in (1/2, 1)$ in (\ref{E:ft}), then
$$
f_t(x) \le \sqrt{x}
$$ 
and the positivity of the matrix
$$
Y_{ij}=m_{f_t}(\lambda_i, \lambda_j)
$$
can be shown similarly to the above argument.

\begin{pl}\label{pl:ando}
The mean
$$
m(x,y)=\frac{1}{2}\left(\frac{x+y}{2}+\frac{2xy}{x+y}\right)
$$ 
is larger than the geometric mean. Indeed,
$$
\frac{1}{2}\left(\frac{x+y}{2}+\frac{2xy}{x+y}\right)\ge
\sqrt{\frac{x+y}{2} \frac{2xy}{x+y}}=\sqrt{xy}.
$$

The numerical computation shows that in this case already the determinant
of a $3 \times 3$ matrix $T$ can be negative. This example shows that the 
corresponding mapping $\beta$ is not completely positive.  \qed
\end{pl}

Next we consider the function
\begin{equation} \label{E:efek}
f_t(x)=t(1-t)\frac{(x-1)^2}{(x^t-1)(x^{1-t}-1)}
\end{equation}
which was first studied in the paper \cite{HP}. If $0 < t < 1$, then the
integral representation
\begin{equation}\label{E:PHas}
{1 \over f_t (x)} ={\sin t \pi \over \pi} \int_0^{\infty} d
\lambda \, \lambda^{t-1}\int_0^1 ds \int_0^1 dr {1 \over x((1-r) \lambda +(1-s))
+(r \lambda +s)}
\end{equation}
shows that  $f_t(x)$ is operator monotone. (Note that in the paper \cite{Szabo} the
operator monotonicity was obtained for $-1 \le t \le 2$.) The property $xf(x^{-1})=
f(x)$ is obvious. 

If $t=1/2$, then 
$$
f(x)=\left(\frac{1+\sqrt{x}}{2}\right)^2 \ge \sqrt{x}
$$
and the corresponding mean is called binomial or power mean. In this case we have
$$
T_{ij}=\frac{4}{(\sqrt{\lambda_i}+\sqrt{\lambda_j})^2}.
$$
The matrix
$$
U_{ij}=\frac{1}{\sqrt{\lambda_i}+\sqrt{\lambda_j}}
$$
is a kind of Cauchy matrix, so it is positive. Since $T=4 U\circ U$, $T$ is positive 
as well. 

If $\gamma (A)=A \circ U$, then $\beta=4 \gamma^2$. Since
$$
\gamma (A)=\int_0^\infty \exp(-s\sqrt{D})A\exp(-s\sqrt{D})\, ds,
$$  
we have
\begin{equation}
\beta(A)=4 \int_0^\infty \int_0^\infty \exp(-(s+r)\sqrt{D})A\exp(-(s+r)\sqrt{D})\, ds\,dr.
\end{equation}
The complete positivity of $\beta$ is clear from this formula.

For the other values of $t$ in $(0,1)$ the proof is a bit more sophisticated.

\begin{lemma}
If $0< t <1$, then $f_t(x)\geq \sqrt{x}$ for $x>0$.
\end{lemma}

\proof
It is enough to show that for $0<t<1$ and $x>0$
\begin{equation}\label{bizfo}
t\frac{x-1}{x^t-1}\geq x^{\frac{1-t}2},
\end{equation}
since this implies
\[
t\frac{x-1}{x^t-1}(1-t)\frac{x-1}{x^{1-t}-1}\geq x^{\frac{1-t}2}x^{\frac{t}2}=\sqrt{x}.
\]

Denote
\[
g(x):=t(x-1)+x^{\frac{1-t}2}-x^{\frac{1+t}2}.
\]
Then inequality (\ref{bizfo}) reduces to $g(x)\geq0$ for $x\geq1$ and to $g(x)\leq0$ 
for $0<x\leq1$. Since $g(1)=0$ it suffices to verify that $g$ is monotone increasing, 
in other words $g'\geq0$. By simple calculation one obtains 
\[g
'(x)=t+\frac{1-t}2x^{\frac{-t-1}2}-\frac{1+t}2x^{\frac{t-1}2}
\]
and
\[
g''(x)=\frac{1-t^2}4x^{\frac{t-3}2}-\frac{1-t^2}4x^{\frac{-t-3}2},
\]
which yields $g''(x)\leq0$ for $0<x<1$ and $g''(x)\geq0$ for $x\geq1$. 
Thus, due to $g'(1)=0$, $g'\geq0$, the statement follows. \qed

It follows from the lemma that the matrix
$$
T_{ij}=t(1-t)\times \frac{\lambda_i^t-\lambda_j^t}{\lambda_i-\lambda_j} \times
\frac{\lambda_i^{1-t}-\lambda_j^{1-t}}{\lambda_i-\lambda_j} \qquad (1 \le i,j \le m) 
$$
can be positive. It is a Hadamard product, so it is enough to see that
$$
U_{ij}^{(t)}=\frac{\lambda_i^t-\lambda_j^t}{\lambda_i-\lambda_j}\qquad (1 \le i,j \le m)  
$$
is positive for $0 < t <1$. It is a well-known fact that the function $g: \bbbr^+ \to 
\bbbr$ is matrix monotone if and only if the L\"owner matrices
$$
L_{ij}=\frac{g(\lambda_i)-g(\lambda_j)}{\lambda_i-\lambda_j}\qquad (1 \le i,j \le m) 
$$ 
are positive. The function $g(x)=x^t$ is matrix monotone for  $0 < t <1$ and the 
positivity of $U$ and $T$ follows. So we have:

\begin{thm}\label{T:2}
For the function (\ref{E:efek}) the mapping $\beta$ is completely positive if $0< t< 1$.
\end{thm}

To see the explicit complete positivity of $\beta$, the mappings $\gamma_t(A)=
A \circ U^{(t)}$ are useful, we have 
$$
\beta (A)= t(1-t)\gamma_t (\gamma_{1-t}(A)). 
$$
Instead of the Hadamard product, which needs the diagonality of $D$, we can use
$$
\gamma_t(A)=\frac{\partial}{\partial x}(D+xA)^t\Big|_{x=0}.
$$
We compute $\gamma_t$ from
$$
(D+xA)^t=\frac{\sin \pi t}{\pi}\int_0^\infty \left(I-s(D+xA+s I)^{-1}\right) 
s^{t-1}\,ds.
$$
So we obtain
$$
\gamma_t(A)=\frac{\sin \pi t}{\pi}\int_0^\infty s^{t}(D+s I)^{-1}A(D+s I)^{-1}
\,ds
$$
and
\begin{eqnarray}
\beta (A)&=&t(1-t)\frac{\sin \pi t  \,\sin \pi (1-t)}{\pi^2}\cr && 
\int_0^\infty \int_0^\infty r^{1-t} s^{t}(D+r I)^{-1}(D+s I)^{-1}A(D+s I)^{-1}(D+r I)^{-1}
\,ds\,dr.
\end{eqnarray}

\begin{pl}\label{Pl:6}
The {\it power difference means} are determined by the functions
\begin{equation}\label{E:ll}
f_t(x)=\frac{t-1}{t}\frac{x^t-1}{x^{t-1}-1} \qquad (-1 \le t \le 2),
\end{equation}
where the values $t=-1,\,1/2,\, 1,\, 2$ correspond to the well-known means, harmonic, 
geometric, logarithmic and arithmetic. The functions (\ref{E:ll}) are operator monotone 
\cite{Fur} and we show that for fixed $x>0$ the value $f_t(x)$ is increasing function
of $t$.

By substituting $x=e^{2\lambda}$ one has
\[
f_t(e^{2\lambda})=\frac{t-1}{t}\frac{e^{\lambda t}\frac{e^{\lambda t}-e^{-\lambda t}}2}{e^{\lambda (t-1)}
\frac{e^{\lambda (t-1)}-e^{-\lambda (t-1)}}2}=e^\lambda\frac{t-1}{t}\frac{\sinh(\lambda t)}
{\sinh(\lambda(t-1))}.
\]
Since 
\[
\frac{d}{dt}\left(\frac{t-1}{t}\frac{\sinh(\lambda t)}{\sinh(\lambda(t-1))}\right)=
\frac{\sinh(\lambda t)\sinh(\lambda(t-1))-\lambda t(t-1)\sinh(\lambda)}
{t^2\sinh^2(\lambda (t-1))},
\]
it suffices to show that
\[
g(t)=\sinh(\lambda t)\sinh(\lambda(t-1))-\lambda t(t-1)\sinh(\lambda)\geq0.
\]
Observe that $\lim_{\pm\infty}g=+\infty$ thus $g$ has a global minimum.
By simple calculations one obtains 
\[
g'(t)=\lambda(\sinh(\lambda(2t-1))-(2t-1)\sinh(\lambda)).\]
It is easily seen that the zeros of $g'$ are $t=0$, $t=1/2$ and $t=1$ hence $g(0)=g(1)=0$ 
and $g(\frac12)=\sinh^2(\frac{\lambda}2)+\frac{\lambda}4\sinh(\lambda)\geq0$ implies 
that $g\geq0$.

It follows that
$$
\sqrt{x} \le f_t(x)\le \frac{1+x}{2}
$$
when $1/2 \le t \le 2$. For these values of the parameter $t$ the complete positivity
holds. This follows from the next lemma which contains a bigger interval for $t$.

\begin{lemma}
The matrix
\[
T_{ij}:=\frac{t}{t-1}\frac{\lambda_i^{t-1}-\lambda_j^{t-1}}{\lambda_i^t-\lambda_j^t}
\]
is positive if $\frac12\leq t$. 
\end{lemma}

\proof
For $t>1$ the statement follows from the proof of Theorem  \ref{T:2}, since
\[
\frac{t}{t-1}\frac{\lambda_i^{t-1}-\lambda_j^{t-1}}{\lambda_i^t-\lambda_j^t}
=\frac{t}{t-1}\frac{(\lambda_i^t)^{\frac{t-1}{t}}-(\lambda_j^t)^{\frac{t-1}{t}}}
{\lambda_i^t-\lambda_j^t},
\]
where $0<\frac{t-1}{t}<1$, further, for $t=1$ the statement follows from Example \ref{P:3}.
If $\frac12\leq t<1$ let $s:=1-t$ where $0<s\leq \frac12$. Then
\[
T_{ij}=\frac{t}{t-1}\frac{\lambda_i^{t-1}-\lambda_j^{t-1}}{\lambda_i^t-\lambda_j^t}
=\frac{1-s}{-s}\frac{\lambda_i^{-s}-\lambda_j^{-s}}{\lambda_i^t-\lambda_j^t}
=\frac{1-s}{s}\frac{(\lambda_i^{t})^{\frac{s}{t}}-(\lambda_j^t)^{\frac{s}{t}}}
{\lambda_i^t-\lambda_j^t}\frac{1}{\lambda_i^{s}\lambda_j^{s}}
\]
so that $T$ is the Hadamard product of $U$ and $V$, where
\[
U_{ij}=\frac{(\lambda_i^t)^{\frac{s}{t}}-(\lambda_j^t)^{\frac{s}{t}}}{\lambda_i^t-\lambda_j^t}
\]
is positive due to $0<\frac{s}{t}\leq1$ and 
\[
V_{ij}=\frac{1-s}{s}\frac{1}{\lambda_i^{s}\lambda_j^{s}}
\]
is positive, too. \qed
\end{pl}

\begin{pl}
Another interpolation between the arithmetic mean  ($t=1$) and the harmonic mean ($t=0$)
is the following:
$$
f_t(x)=\frac{2(tx+1)(t+x)}{(1+t)^2(x+1)} \qquad (0 \le t \le 1).
$$

First we compare this mean with the geometric mean:
$$
f_t(x^2)-x=\frac{(x-1)^2(2tx^2-(1-t)^2x +2t)}{(1+t)^2(x^2+1)}
$$
and the sign depends on 
$$
x^2-\frac{(1-t)^2}{2t}x + 1=\left(x-\frac{(1-t)^2}{4t}\right)^2+1-
\left(\frac{(1-t)^2}{4t}\right)^2.
$$
So the positivity condition is $(1-t)^2 \le 4t$ which gives $3-2 \sqrt{2}\le t
\le 3+2 \sqrt{2}$.
For these parameters $f_t(x) \ge \sqrt{x}$ and for $0<t < 3-2 \sqrt{2}$ the two
means are not comparable.

For $3-2 \sqrt{2}\le t \le 1$ the matrix monotonicity is rather straightforward:
$$
f_t(x)=\frac{2}{(1+t)^2}\left(tx+t^2-t+1-\frac{(t-1)^2}{x+1}\right)
$$
However, the numerical computations show that $T \ge 0$ is not true. \qed
\end{pl}

%%%%%%%%%%%%%%%%%%%%%%%%%%%%%%%%%%%%%%%
\section{Some matrix monotone functions}

First the {\it Stolarsky mean} is investigated \cite{Pales, Sto}.

\begin{thm}
Let
\begin{equation} 
f_p(x):=\left(\frac{p(x-1)}{x^p-1}\right)^{\frac1{1-p}}, 
\end{equation}
where $p \ne 1$. Then $f_p$ is matrix monotone if $-2\leq p\leq 2$.
\end{thm}

\proof
First note that $f_2(x)=(x+1)/2$ is the arithmetic mean, the limiting case 
$f_0(x)=(x-1)/\log x$ is the logarithmic mean and $f_{-1}(x)=\sqrt{x}$ is 
the geometric mean, their matrix monotonicity is well-known. If $p=-2$ then 
$$
f_{-2}(x)=\frac{(2x)^{\frac23}}{(x+1)^{\frac13}}
$$ 
which will be shown to be matrix monotone at the end of the proof.

Now let us suppose that $p\neq -2, -1,0,1,2$. By L\"owner's theorem $f_p$ is matrix 
monotone if and only if it has a holomorphic continuation mapping the upper half 
plane into itself. We define $\log z$ as $\log1:=0$ then in case $-2<p<2$, since 
$z^p-1\neq0$ in the upper half plane, the real function $p(x-1)/(x^p-1)$ 
has a holomorphic continuation to the upper half plane, moreover it is continuous 
in the closed upper half plane, further, $p(z-1)/(z^p-1)\neq0$ 
($z\neq1$) so $f_p$ also has a holomorphic continuation to the upper half plane 
and it is also continuous in the closed upper half plane. 

Assume $-2<p<2$ then it suffices to show that $f_p$ maps the upper half plane into itself. 
We show that for every $\varepsilon>0$ there is $R>0$ such that the set $\{z:|z|\geq R,
\Imag z>0\}$ is mapped into $\{z: 0\leq \arg z\leq\pi+\varepsilon\}$, further, the 
boundary $(-\infty,+\infty)$ is mapped into the closed upper half plane. Then by the 
well-known fact that the image of a connected open set by a holomorphic function is 
either a connected open set or a single point it follows that the upper half plane 
is mapped into itself by $f_p$.

Clearly, $[0,+\infty)$ is mapped into $[0,\infty)$ by $f_p$. 

Now first suppose $0<p<2$. Let $\varepsilon>0$ be sufficiently small and 
$z\in\{z:|z|=R,\,\Imag z>0\}$ where $R>0$ is sufficiently large. Then
\[
\arg (z^p-1)=\arg z^p\pm\varepsilon=p\arg z\pm \varepsilon,
\]
and similarly $\arg z-1=\arg z\pm\varepsilon$ so that 
\[
\arg \frac{z-1}{z^p-1}=(1-p)\arg z\pm 2\varepsilon.
\]
Further,
\[
\left|\frac{z-1}{z^p-1}\right|\geq \frac{|z|-1}{|z|^p+1}=\frac{R-1}{R^p+1},
\]
which is large for $0<p<1$ and small for $1<p<2$ if $R$ is sufficiently large, hence
\[
\arg\left(\frac{z-1}{z^p-1}\right)^{\frac1{1-p}}=\frac1{1-p}\arg
\left(\frac{z-1}{z^p-1}\right)\pm2\varepsilon=\arg z\pm2\varepsilon\frac{2-p}{1-p}.
\]
Since $\varepsilon>0$ was arbitrary it follows that $\{z:|z|=R,\,\Imag z>0\}$ is mapped 
into the upper half plane by $f_p$ if $R>0$ is sufficiently large.

Now, if $z\in [-R,0)$ then $\arg (z-1)=\pi$, further, $p\pi\leq \arg (z^p-1)
\leq\pi$ for $0<p<1$ and $\pi\leq \arg (z^p-1)\leq p\pi$ for $1<p<2$ whence 
\[
0\leq\arg\left(\frac{z-1}{z^p-1}\right)\leq (1-p)\pi\quad\mbox{for}\quad 0<p<1,
\] 
and
\[
(1-p)\pi\leq\arg\left(\frac{z-1}{z^p-1}\right)\leq 0\quad\mbox{for}\quad 1<p<2.
\] 
Thus by
\[
\pi\arg\left(\frac{z-1}{z^p-1}\right)^{\frac1{1-p}}=\frac1{1-p}\arg\left(\frac{z-1}{z^p-1}\right)
\]
it follows that
\[
0\leq\arg\left(\frac{z-1}{z^p-1}\right)^{\frac1{1-p}}\leq\pi
\]
so $z$ is mapped into the closed upper half plane.

The case $-2<p<0$ can be treated similarly by studying the arguments and noting that
\[
f_p(x)=\left(\frac{p(x-1)}{x^p-1}\right)^{\frac1{1-p}}=\left(\frac{|p|x^{|p|}(x-1)}{x^{|p|}-1}
\right)^{\frac1{1+|p|}}.
\]

Finally, we show that $f_{-2}(x)$ is matrix monotone. 
Clearly $f_{-2}$ has a holomorphic continuation to the upper half plane (which is not continuous in the closed upper half plane). If $0<\arg z<\pi$ then $\arg z^{\frac23}=\frac23\arg z$ and $0<\arg 
(z+1)<\arg z$ so 
\[
0<\arg\left(\frac{z^{\frac23}}{(z+1)^{\frac13}}\right)<\pi
\]
thus the upper half plane is mapped into itself by $f_{-2}$. \qed

The limiting case $p=1$ is the so-called {\it identric mean}:
\[
f_1(x)=\frac1e x^{\frac{x}{x-1}}=\exp\left(\frac{x\log x}{x-1}-1\right).
\]
It is not so difficult to show that $f_1$ is matrix monotone.

The inequality
$$
\sqrt{x} \le f_p(x) \le \frac{1+x}{2}
$$
holds if $p \in [-1, 2]$. It is proved in \cite{Kos} that the matrix
$$
T_{ij}=\left(\frac{\lambda_i^p-\lambda_j^p}{p(\lambda_i-\lambda_j)}\right)^{\frac1{1-p}}
$$
is positive. 

\begin{cor}
The mapping $\beta$ induced by the Stolarsky mean is monotone and  completely positive
for $p \in [-1, 2]$.  
\end{cor}

The {\it power} or {\it binomial mean}
$$
m(a,b)=\left(\frac{a^p+b^p}2\right)^{\frac1p}
$$
can be also a matrix monotone function:

\begin{thm}
The function
\begin{equation} 
f_p(x)=\left(\frac{x^p+1}2\right)^{\frac1p}
\end{equation}
is matrix monotone if and only if $-1\leq p\leq 1$.
\end{thm}

\proof
Observe that $f_{-1}(x)=2x/(x+1)$ and $f_1(x)=(x+1)/2$, so $f_p$ could be 
matrix monotone only if $-1\leq p\leq1$. We show that it is indeed matrix monotone. 
The case $p=0$ is well-known. Further, note that if $f_p$ is matrix monotone for $0<p<1$ then 
\[
f_{-p}(x)=\left(\left(\frac{x^{-p}+1}2\right)^{\frac1p}\right)^{-1}
\]
is also matrix monotone since $x^{-p}$ is matrix monotone decreasing for $0<p\leq1$.

So let us assume that $0<p<1$. Then, since $z^p+1\neq0$ in the upper half plane, 
$f_p$ has a holomorphic continuation to the upper half plane (by defining $\log z$ as $\log 1=0$). By L\"owner's theorem 
it suffices to show that $f_p$ maps the upper half plane into itself. If $0<\arg z<
\pi$ then $0<\arg (z^p+1)<\arg z^p=p\arg z$ so
\[
0<\arg \left(\frac{z^p+1}2\right)^{\frac1p}=\frac1p\arg\left(\frac{z^p+1}2\right)<\arg z<\pi
\]
thus $z$ is mapped into the upper half plane. \qed

In the special case $p=\frac1n$,
\[
f_p(x)=\left(\frac{x^{\frac1n}+1}2\right)^n=\frac1{2^n}\sum_{k=0}^n{n \choose k}x^{\frac{k}n},
\]
and it is well-known that $x^\alpha$ is matrix monotone for $0<\alpha<1$ thus $f_p$ is 
also matrix monotone.

Since the power mean is infinitely divisible \cite{Kos}, we have:

\begin{cor}
The mapping $\beta$ induced by the power mean is monotone and  completely positive
for $p \in [-1, 1]$.  
\end{cor}

\section{Discussion and conclusion}

The complete positivity of some mappings $\beta_f : \bM_n \to \bM_n$ has been
a question in physical applications when $\beta$ is determined by a function
$f: \bbbr^+ \to  \bbbr^+$. The function $f$ is in connection with means
larger than the geometric mean. In the paper several concrete functions
are studied, for example, Heinz mean, power difference means, Stolarsky mean 
and interpolations between some means. The complete positivity of $\beta_f$ 
is equivalent with the positivity of a matrix. The analysis of the functions 
studied here is very concrete, general statement is not known. 
 
{\bf Acknowledgements.} Thanks to several colleagues for communication. 
Professor Tsuyoshi Ando suggested Example \ref{pl:ando}, the student  G\'abor 
Ball\'o made numerical computations, Professor Zsolt P\'ales explained  
certain means and Professor M. Beth Ruskai informed about the question of \cite{TRus}.


\begin{thebibliography}{99}

\bibitem{Kos}
R. Bhatia and H. Kosaki, Mean matrices and infinite divisibility,  
Linear Algebra Appl.  {\bf 424}(2007), 36--54.

\bibitem{Fur}
T. Furuta, Concrete examples of operator monotone functions obtained by an elementary 
method without appealing to L\"owner integral representation,  Linear Algebra Appl.  
{\bf 429}(2008), 972--980.

\bibitem{Han}
F. Hansen, Characterization of symmetric monotone metrics on the state space of 
quantum systems, Quantum Inf. Comput.  {\bf 6}(2006), 597--605. 

\bibitem{HP}
H. Hasegawa and D. Petz,
On the Riemannian metric of $\alpha$-entropies of density matrices,
Lett. Math. Phys. {\bf 38}(1996), 221--225.

\bibitem{HK}
F. Hiai and H. Kosaki, Comparison of various means for operators,  
J. Funct. Anal.  {\bf 163}(1999), 300--323.

\bibitem{HK2}
F. Hiai and H. Kosaki, {\it Means of Hilbert space operators}, Lecture Notes 
in Mathematics, 1820. Springer-Verlag, Berlin, 2003. 

\bibitem{PD145}
F. Hiai and D. Petz, From quasi-entropy, arXiv:1009.2679, 2010.

\bibitem{K-A}
F. Kubo and T. Ando, Means of positive linear operators, Math. Ann. {\bf 246}(1980), 
205--224.

\bibitem{Pales}
E. Neumann and Zs. P\'ales, On comparison of Stolarsky and Gini means, J. Math. Anal. Appl.,
{\bf 278}(2003), 274--285.

\bibitem{PD2}
D. Petz, Monotone metrics on matrix spaces, Linear
Algebra Appl. {\bf 244}(1996), 81--96.

\bibitem{PD22}
D. Petz,  Covariance and Fisher information in quantum mechanics,
J. Phys. A: Math. Gen. {\bf 35}(2003), 79--91.

\bibitem{PD23}
D. Petz, Means of positive matrices: Geometry and a conjecture,
Annales Mathematicae et Informaticae {\bf 32}(2005), 129--139.

\bibitem{PD143}
D. Petz and C. Ghinea, Introduction to quantum Fisher information, arXiv:1008.2417,
to appear in {\it QP--PQ: Quantum Probab. White Noise Anal., vol. 27.}

\bibitem{Sto}
K.B. Stolarsky, Generalizations of the logarithmic mean,  Math. Mag.  
{\bf 48}(1975), 87--92.

\bibitem{Szabo}
V.E.S. Szab\'o, A class of matrix monotone functions,  Linear Algebra Appl. 
{\bf 420}(2007), 79--85. 

\bibitem{TRus}
K. Temme, M. J. Kastoryano, M. B. Ruskai, M. M. Wolf and F. Verstraete,
The $\chi^2$-divergence and mixing times of quantum Markov processes,
J. Math. Phys. {\bf 51}(2010), 122201.  


\end{thebibliography}
\end{document}